\documentclass{smathjsty}
\newtheorem{theorem}{Theorem}[section]

\newtheorem{proposition}[theorem]{Proposition}

\theoremstyle{definition}
\newtheorem{definition}[theorem]{Definition}
\newtheorem{example}[theorem]{Example}

\theoremstyle{remark}
\newtheorem{remark}[theorem]{Remark}
\theoremstyle{proof}
\numberwithin{equation}{section}

\begin{document}

\DateofSubmission{May 21, 2011}
\EmailofCorrespondingAuthor{amah@vru.ac.ir; mohsenialhosseini@gmail.com}

\title[G-approximate best proximity pairs
in metric space with a directed graph] {G-approximate best proximity pairs
in metric space with a directed graph}
\author[S. A. M. Mohsenailhosseini, M. Saheli] {S. A. M. Mohsenialhosseini$^*$, M. Saheli}
\maketitle
\vspace*{-0.2cm}
\begin{center}
{\footnotesize  
 $^*$
 Faculty of Mathematics, Vali-e-Asr University, Rafsanjan, Iran. e-mail: amah@vru.ac.ir\\
 Faculty of Mathematics, Vali-e-Asr University, Rafsanjan, Iran. e-mail: saheli@vru.ac.ir
}
\end{center}

{\tiny\CopyrightInfo}

{\footnotesize 

\noindent {\bf Abstract.} 
Let $(X,d)$ be a metric space that has a directed graph G such that the sets $V(G)$ and $E(G)$ are respectively vertices and edges corresponding to $X.$ We obtain sufficient conditions for the existence of an $G$-approximate best proximity pair of the mapping $T$ in the metric space $X$ endowed with a graph $G$ such that the set $V(G)$ of vertices of $G$ coincides with $X.$
\vskip 1mm


}

\vskip 6mm

\FootNoteInfo 
\section { Introduction  }
 Fixed point theory is a very popular tool in solving existence
problems in many branches of Mathematical Analysis and its applications. In physics and engineering fixed point technique has been used in areas like image retrieval,  signal processing and the study of existence and uniqueness of solutions for a class of nonlinear integral equations.\\
Graphs can be used to model many types of relations and processes in physical, biochemistry, electrical engineering, computer science  and operations research, biological.The wide scope of these and other applications has been well-documented cf. \cite {Cac,Rob}
\par Let X be a metric space and A and B nonempty subsets of X and $d(A,B)$ is distance of $A$ and $B$. If  $d(x_0,y_0)=d(A,B)$ then the pair $(x_0,y_0)$ is called a best proximity pair for A and B and put $$prox(A,B):=\{(x,y)\in A\times B :~d(x,y)=d(A,B)\}$$ as the set of all best proximity pair (A,B). Best proximity pair evolves as a generalization of the concept of best approximation. That reader can find some important result of it in \cite {Fan,kim,kir}. \vspace{0.4cm}

\par Now, as in \cite {Sin} (see also\cite {Eld,Bee,Ve,ple,ple1,ple2}), we can find the best proximity points of the sets A and B, by considering a map   $T:A\cup B \rightarrow A\cup B$ such that $T(A)\subseteq B$ and  $T(B)\subseteq A$. Best proximity pair also evolves  as a generalization of the concept of fixed point of mappings. Because if $A\cap B\neq \emptyset$ every best proximity point is a fixed point of $T$.\\
In 2011, Mohsenalhosseini et al \cite {Moh0},  introduced to approximate best proximity pairs in metric space.  In 2013, Mohsenalhosseini et al \cite {Moh1},  introduced  to approximate best proximity pairs  in metric space for contraction maps.\\
Recently, two results have appeared, giving suficient conditions for f to be a Picard operator
if $(X, d)$ is endowed with a graph. The first result in this direction was given by J. Jakhymski \cite {Jac} who also presented its applications to the Kelisky-Rivlin theorem on iterates of the Bernstein operators on the space $C [0,1]$. \\
Reich \cite {Rei}, \'Ciri\'c \cite {Cir} and Rus \cite {Rus1} proved that if (X, d) is complete, then every \'Ciri\'c– Reich–Rus operator has a unique fixed point. 

\par In \cite {Boj}, the author considered the problem of existence of a fixed point for $\varphi$-contractions in metric spaces endowed with a graph. Also in \cite {Boj1}, introduced to Fixed point theorems for Reich type contractions on metric spaces with a graph.
\par The aim of this paper is to study the existence of of approximate best proximity pair for a map and two maps and their diameters for a cyclic map  $T:A\cup B\cup C\rightarrow A\cup B\cup C$ i.e. $T(A)\subseteq B$, $T(B)\subseteq C$ and $T(C)\subseteq A$ in metric spaces endowed with a graph $G$ by defining the G-approximate best proximity pair. Moreover, we give some illustrative example of our main results. 

\section { Preliminaries }
This section recalls the following notations and the ones that will be used in what follows.
We refer to Mohsenalhosseini et al \cite {Moh0, Moh1}  for a detailed study of approximate best proximity pair for a map and two maps,  also refer to Wayne Goddard et al \cite {God} for a detailed study of Distance in Graphset.\\
Graphs are special examples of metric spaces with their intrinsic path metric. 
Now we express the relationship between distance, diameter and radius in the graph with the metric space.
A path in a graph is a sequence of distinct vertices, such that adjacent vertices in the sequence are adjacent in the graph. For an unweighted graph, the length of a path is the number of edges on the path. For an (edge) weighted graph, the length of a path is the sum of the weights of the edges on the path. We assume that all weights are nonnegative and that all graphs are connected. 

\begin{definition} \cite {Moh2} 
{\it Let $T:X\rightarrow X,$ $\epsilon>0,$  $x_0\in X.$
 Then $x_0\in X$ is an $\epsilon-$fixed point for $T$ if $d(x_0,Tx_0)<\epsilon.$}
\end{definition}
\begin{definition} \cite {Moh2} Let $T:X\rightarrow X.$  Then $T$ has the approximate fixed point property (a.f.p.p) if $$\forall \epsilon>0,~F_\epsilon (T)\neq\varnothing.$$
\end{definition}

\begin{theorem} \cite {Moh1} \label{1.6} {\it Let $(X,\|.\|)$ be a  complete norm space,
$T:X\rightarrow X,$   $x_0\in X$ and $\epsilon>0$ . If $\|T^{n}(x_0)-T^{n+k}(x_0)\|\rightarrow 0$ as $n\rightarrow\infty$ for some $k>0,$  then $T^k$ has an $\epsilon-$ fixed point.}
\end{theorem}

\begin{definition} \cite {Jac}\label{1.1} We say that a mapping  $T:X\rightarrow X$ is a $G$-contraction or
simply $G$-contraction if $T$ preserves edges of $G$, i.e.,
\begin{equation}
\forall x,y\in X ((x,y)\in E(G)\Rightarrow (Tx,Ty)\in E(G)),
\end{equation}
and $T$ decreases weights of edges of $G$ in the following way:
\begin{equation}
\exists \alpha \in(0,1)~ \forall x,y\in X~ ((x,y)\in E(G)\Rightarrow d(Tx,Ty)\leq \alpha d(x,y)).
\end{equation}
\end{definition}

\section { Main results }

Throughout this section, we assume that  $G$ is a directed graph such that $V(G) = X$, $E(G) \supseteq \bigtriangleup$ such that  $\bigtriangleup$ denote the diagonal of the Cartesian product $X \times X$ and the graph $G$ has no parallel edges. In this section, by using the idea of Jachymski \cite {Jac}, we will consider the existence of $G$-approximate best proximity points for the map $T:~A \cup B\rightarrow A \cup B$, such that  $T(A)\subseteq B$, $T(B)\subseteq A$, and its diameter.
\begin{definition} \label{2.1}Let $A$ and $B$ be nonempty subsets of a metric space endowed with a directed graph $G=(V(G),E(G))$ such that $V(G)=A\cup B$. The operator $T:A\cup B \rightarrow A\cup B$ be a map such that  $T(A)\subseteq B$, $T(B)\subseteq A.$  The point $x\in A\cup B$ is said to be a $G$-approximate best proximity point of the pair $(A,B)$, if:\\
 i) $((x,y)\in E(G)\Rightarrow (Tx,Ty)\in E(G)), \forall x,y \in A\cup B$;\\
ii) $\exists\epsilon>0~ \forall x\in A\cup B~ ((x,Tx)\in E(G)\Rightarrow d(x,Tx)\leq d(A,B)+\epsilon).$ 
\end{definition}
\begin{remark}   \label{1.2}
In this paper we will denote the set of all $G$-approximate best proximity point of the pair $(A,B)$ of $T$, for a
given $\epsilon$, by : 
\begin{eqnarray*}
P^{Ga}_{T}(A,B)=\{ x\in A\cup B:((x,Tx)\in E(G)\Rightarrow d(x,Tx)\leq d(A,B)+\epsilon)\}.
\end{eqnarray*}
\end{remark}
\begin{example}
Suppose Let $X=\textbf{R}^2$ and $A=\{(x,y)\in X:(x-y)^2+y^2\leq 1\}$ and $B=\{(x,y)\in X:(x+y)^2+y^2\leq 1\}$ with $T(x,y)=(-x,y)$ for $(x,y)\in X$. Define the graph $G$ by $V(G)=X.$
Then there exists $((x,y),T(x,y))\in E(G)$ such that  $d((x,y),T(x,y))\leq d(A,B)+\epsilon$  for some $\epsilon>0$. Hence $P^{Ga}_{T}(A,B)\neq\emptyset.$ 
\end{example}
We say that the pair $(A,B)$ is a $G$-approximate best proximity pair if $P^{Ga}_{T}(A,B)\neq\emptyset.$

\begin{proposition} \label{3.4} Let $A$ and $B$ be nonempty subsets of a metric space endowed with a directed graph $G=(V(G),E(G))$ such that $V(G)=A\cup B$. suppose that the mapping $T:A\cup B \rightarrow A\cup B$ satisfying   $T(A)\subseteq B$, $T(B)\subseteq A.$ If $\lim_{n\rightarrow\infty}d(T^nx,T^{n+1}x)=d(A,B)$, for some $x \in A\cup B$ satisfies the condition $(x, Tx)\in E(G)$ then the pair $(A,B)$ is a  $G$-approximate best proximity pair.
\end{proposition}
{\bf Proof:} Let  $\epsilon>0$ be given and  $x \in A\cup B$ with $(x, Tx)\in E(G)$ such that $\lim_{n\rightarrow\infty}d(T^nx,T^{n+1}x)=d(A,B)$; then there exists $N_{0}>0$ such that for all $n\geq N_0,$ \[d(T^nx,T^{n+1}x)<d(A,B)+\epsilon.\]
If $n=N_0,$ then $d(T^{N_0}(x),T(T^{N_0}(x)))< d(A,B)+\epsilon,$
then $T^{N_0}(x)\in P^{Ga}_{T}(A,B)$ and $P^{Ga}_{T}(A,B)\neq\emptyset$. $\blacksquare $
\begin{definition}  Let $A$ and $B$ be nonempty subsets of a metric space endowed with a graph $G$. The operator $T:A\cup B \rightarrow A\cup B$ satisfying   $T(A)\subseteq B$, $T(B)\subseteq A$  is said to be a G-ciric-Rich-Rus-Moh operator if:\\
 i) $ ((x,y)\in E(G)\Rightarrow (Tx,Ty)\in E(G)), \forall x,y \in A\cup B$;\\
ii) there exists nonnegative numbers $\alpha,\beta,\gamma$ with $\alpha+2\beta+\gamma<1$, such that, for each $(x, y)\in E(G),$
 we have:$$d(Tx,Ty)\leq \alpha d(x,y)+\beta[d(x,Tx)+d(y,Ty)]+\gamma d(A,B)$$
\end{definition}
\begin{theorem}   \label{3.9}
 Let $A$ and $B$ be nonempty subsets of a metric space endowed with a graph $G$ and $T:A\cup B \rightarrow A\cup B$ be a  G-ciric-Rich-Rus-Moh operator. If $x \in A\cup B$ satisfies the condition $(x, Tx)\in E(G)$ then the pair $(A,B)$ is a $G$-approximate best proximity pair.
 \end{theorem}
 
 {\bf Proof:} Let  $x \in A\cup B$ with $(x, Tx)\in E(G)$, then $$d(Tx,T^2x)\leq \alpha d(x,Tx)+\beta[d(x,Tx)+d(Tx,T^2x)]+\gamma d(A,B)$$ Therefore $$d(Tx,T^2x)\leq \frac{\alpha+\beta}{1-\beta} d(x,Tx)+\frac{\gamma}{1-\beta} d(A,B).$$
Now if $k=\frac{\alpha+\beta}{1-\beta},$ then $$d(Tx,T^2x)\leq kd(x,Tx)+(1-k)d(A,B)$$ also $$d(T^2x,T^3x)\leq k^2d(x,Tx)+(1-k^2)d(A,B).$$ Therefore $$d(T^nx,T^{n+1}x)\leq k^nd(x,Tx)+(1-k^n)d(A,B),$$ and so $$d(T^nx,T^{n+1}x)\rightarrow d(A,B),~as~ n\rightarrow\infty.$$
Therefore, by Theorem \ref{3.4}, $P^{Ga}_{T}(A,B)\neq\emptyset$, then pair$(A,B)$ is a $G$-approximate best proximity pair. $\blacksquare $
\begin{definition} 
 Let $A$ and $B$ be nonempty subsets of a metric space endowed with a graph $G$. Suppose that the mapping  $T:A\cup B \rightarrow A\cup B$ satisfying   $T(A)\subseteq B$, $T(B)\subseteq A.$ We say that the sequence $\{z_n\}\subseteq A\cup B$ with $(z_n,Tz_n)\in E(G)$ is $GT$-minimizing if $$\lim_{n\rightarrow\infty}d(z_n,Tz_n)=d(A,B).$$
\end{definition}
\begin{theorem} Let $A$ and $B$ be nonempty subsets of a metric space endowed with a graph $G$, suppose that the mapping  $T:A\cup B \rightarrow A\cup B$ satisfying   $T(A)\subseteq B$, $T(B)\subseteq A$. If $\{T^nx\}$ is a $GT$-minimizing  for some $x \in A\cup B$ satisfies the condition $(x, Tx)\in E(G)$ , then $(A,B)$ is a $G$-approximate best pair proximity.
 \end{theorem}
 
{\bf Proof.} Since $$\lim_{n\rightarrow\infty}d(T^nx,T^{n+1}x)=d(A,B)$$ for some $x \in A\cup B$ satisfies the condition $(x, Tx)\in E(G),$ then by Theorem \ref{3.4} $P^{Ga}_{T}(A,B)\neq\emptyset$. Therefore pair $(A,B)$ is a $G$-approximate best proximity pair. $\blacksquare $

\begin{theorem} Let $A$ and $B$ be nonempty subsets of a metric space endowed with a graph $G$ such that $E(G)$ is compact. Suppose that the mapping  $T:A\cup B \rightarrow A\cup B$ satisfying   $T(A)\subseteq B$, $T(B)\subseteq A$, $T$ is continuous and $\|Tx-Ty\|\leq \ \|x-y\|,$ where $(x,y)\in E(G)$. Then $P^{Ga}_{T}(A,B)$ is nonempty and compact.
\end{theorem}

{\bf Proof.} Since $E(G)$ compact, there exists a $z_0\in E(G)$ such that
\begin{equation}
\|z_0-Tz_0\|=\inf_{z\in E(G)}\|z-Tz\|.\label{(3.1)}
\end{equation}
 If $\|z_0-Tz_0\|>d(A,B)$
, then $\|Tz_0-T^2z_0\|<\|z_0-Tz_0\|$ which  contradict to the definition of $z_0$, ($Tz_0\in E(G))$ and by 
\ref{(3.1)}. Therefore $\|z_0-Tz_0\|=d(A,B)\leq d(A,B)+\epsilon$ for some $\epsilon >0$ and $z_0\in P^{Ga}_{T}(A,B)$ . Therefore  $P^{Ga}_{T}(A,B)$ is nonempty.\\
 Also, if $\{z_n\}\subseteq P^{Ga}_{T}(A,B),$ then  $\|z_n-Tz_n\|<d(A,B)+\epsilon,$ for some $\epsilon >0$, and by compactness of $E(G)$, there exists a subsequence ${z_{n_k}}$ and a $z_0\in E(G)$ such that $z_{n_k}\rightarrow z_0$ and so
$$ \|z_0-Tz_0\|=\lim_{k\rightarrow\infty}\|z_{n_k}-Tz_{n_k}\|<d(A,B)+\epsilon$$ for some $\epsilon >0$ hence $P^{Ga}_{T}(A,B)$ is compact. $\blacksquare $

\begin{example} If $A=[-3,-1], B=[1,3]$ and $T:A\cup B \rightarrow A\cup B$ such that
\[T(x)=\left\{\begin {array}{cl}
\frac{1-x}{2}& x\in A \\\\ \frac{-1-x}{2}&x\in B
\end{array}\right.\]
Then $P^{Ga}_{T}(A,B)$ is compact, we have
\begin{eqnarray*}
P^{Ga}_{T}(A,B)&=&\{x\in A\cup B:~d(x,Tx)<d(A,B)+\epsilon~~for~~ some~~ \epsilon>0\}\\
&=& \{x\in A\cup B:~d(x,Tx)<2+\epsilon~~for~~ some~~ \epsilon>0\}\\
&=& \{1,-1\}.
\end{eqnarray*}
that is compact.
\end{example}

\par In the following, by $diam(P^{Ga}_{T}(A,B))$ for a set $P^{Ga}_{T}(A,B)\neq\emptyset$ we will understand the diameter of the set $P^{Ga}_{T}(A,B)$.
\begin{definition} 
  Let $A$ and $B$ be nonempty subsets of a metric space endowed with a graph $G$ such that  $T(A)\subseteq B$, $T(B)\subseteq A$ and $\epsilon >0$. We define diameter $P^{Ga}_{T}(A,B)$ by
 \begin{equation}
 diam(P^{Ga}_{T}(A,B))=\sup \{d(x,y):~~x,y\in P^{Ga}_{T}(A,B)\}.
 \end{equation}
\end{definition} 
\begin{theorem}  
 Let $A$ and $B$ be nonempty subsets of a metric space endowed with a graph $G.$ suppose that the mapping $T:A\cup B \rightarrow A\cup B$, such that  $T(A)\subseteq B$, $T(B)\subseteq A$ and $\epsilon >0$. If $T$ be a $G$-contraction then 
 \begin{equation*}
diam (P^{Ga}_{T}(A,B)) \leq \frac{2\epsilon}{1-\alpha}+\frac{2d(A,B)}{1-\alpha}.
 \end{equation*}
 \end{theorem}

 {\bf Proof.} If $x, y\in P^{Ga}_{T}(A,B)$, then
\begin{eqnarray*}
d(x, y) &\leq& d(x, Tx)+d(Tx, Ty)+d(Ty, y)\\
  &\leq& \epsilon_1 +\alpha d(x, y)+2d(A,B)+\epsilon_2.
\end{eqnarray*}
 put $\epsilon=Max\{\epsilon_1,\epsilon_2\}$, Therefore
 $d(x, y) \leq \frac{2\epsilon}{1-\alpha}+\frac{2d(A,B)}{1-\alpha}.$
Hence $diam (P^{Ga}_{T}(A,B)) \leq \frac{2\epsilon}{1-\alpha}+\frac{2d(A,B)}{1-\alpha}.~\blacksquare $
{\bf 3. {$G$-approximate Best proximity for two maps} }\vspace{0.4cm}

\par In this section we will consider the existence of $G$-approximate best proximity points for two maps $T:A \cup B\rightarrow A \cup B$, $S:A\cup B \rightarrow A\cup B$, and its diameter.
\begin{definition} 
  Let $A$ and $B$ be nonempty subsets of a metric space endowed with a graph $G$ and $T:A\cup B \rightarrow A\cup B$, $S:A\cup B \rightarrow A\cup B$ be two maps such that  $T(A)\subseteq B$, $S(B)\subseteq A$. A point $(x,y)$ in $A\times B$ is said to be a $G$-approximate-pair fixed point for $(T,S)$, if:\\ 
i) $\forall (x,y) \in A \times B~~ ((x,y)\in E(G)\Rightarrow (Tx,Ty)\in E(G), (Sx,Sy)\in E(G))$;\\
ii) $\exists \alpha \in(0,1)~ \forall (x,y)\in A\times B~ ((x,y)\in E(G)\Rightarrow d(Tx,Sy)\leq d(A,B)+\epsilon).$
\end{definition}  
  
We say that the pair $(T,S)$ has the $G$-approximate-pair fixed property in X if   
  $P^{Ga}_{(T,S)}(A,B)\neq \emptyset$, where $$P^{Ga}_{(T,S)}(A,B)=\{(x,y)\in E(G):~d(Tx,Sy)\leq d(A,B)+\epsilon ~~for~ some ~ \epsilon>0\}.$$
 
\begin{theorem} 
 Let $A$ and $B$ be nonempty subsets of a metric space endowed with a graph $G$ and $T:A\cup B \rightarrow A\cup B$, $S:A\cup B \rightarrow A\cup B$ be two maps such that  $T(A)\subseteq B$, $S(B)\subseteq A.$ If, for every $(x,y)\in A\times B$, \[d(T^n(x),S^n(y))
\rightarrow d(A,B),\] then $(T, S)$ has the $G$-approximate-pair fixed property.
\end{theorem}

{\bf Proof.}  Let  $\epsilon>0$ be given and   $(x,y)\in A\times B$ with $((x,y), (Tx,Ty))\in E(G)$. Since $d(T^n(x),S^n(y))
\rightarrow d(A,B),$
there existe $n_0>0$ such that for all  $n\geq n_0,$ \[d(T^n(x),S^n(y)) <d(A,B)+\epsilon.\]
Then $d(T(T^{n-1}(x),S(S^{n-1}(y))<d(A,B)+\epsilon$ for every $n\geq n_0$.
Put $x_0=T^{n_0-1}(x)$ and $y_0=S^{n_0-1}(y))$. Hence $d(T(x_0), S(y_0))\leq d(A,B)+\epsilon$ and $P^{Ga}_{(T,S)}(A,B)\neq \emptyset.~\blacksquare$
\begin{definition}  Let $A$ and $B$ be nonempty subsets of a metric space endowed with a graph $G$. The operator $T:A\cup B \rightarrow A\cup B$, $S:A\cup B \rightarrow A\cup B$ satisfying   $T(A)\subseteq B$, $S(B)\subseteq A$  is said to be a G-ciric-Rich-Rus-2map operator if:\\
 i) $\forall (x,y) \in A \times B~~ ((x,y)\in E(G)\Rightarrow (Tx,Ty)\in E(G), (Sx,Sy)\in E(G))$;\\
ii) there exists nonnegative numbers $\alpha,\beta,\gamma$ with $\alpha+2\beta+\gamma<1$, such that, for each $(x, y)\in E(G),$  we have:$$d(Tx,Sy)\leq \alpha d(x,y)+\beta[d(x,Tx)+d(y,Sy)]+\gamma d(A,B).$$
\end{definition}

\begin{theorem}
 Let $A$ and $B$ be nonempty subsets of a metric space endowed with a graph $G$ and $T:A\cup B \rightarrow A\cup B$, $S:A\cup B \rightarrow A\cup B$ are a  G-ciric-Rich-Rus-2map operator. If $x,y \in A\cup B$ satisfies the condition $(x, Tx),(x, Sx)\in E(G)$, $(y, Ty),(y, Sy)\in E(G)$ and if $x$ is a $G$-approximate fixed point for $T$, or $y$ is a $G$-approximate fixed point for $S$, then $P^{Ga}_{(T,S)}(A,B)\neq \emptyset.$
\end{theorem}

 {\bf Proof.}  Let  $x,y \in A\cup B$ with $(x, Tx),(x, Sx)\in E(G)$, $(y, Ty),(y, Sy)\in E(G)$, then $$d(Tx,S(Tx))\leq \alpha d(x,Tx)+\beta[d(x,Tx)+d(Tx,S(Tx)]+\gamma d(A,B)$$ Therefore $$d(Tx,S(Tx))\leq \frac{\alpha+\beta}{1-\beta} d(x,Tx)+\frac{\gamma}{1-\beta} d(A,B).$$
Now if $k=\frac{\alpha+\beta}{1-\beta},$ then
\begin{equation}\label{3.3}
d(Tx,S(Tx))\leq kd(x,Tx)+(1-k)d(A,B)
 \end{equation}
\begin{equation}\label{3.4}
d(Sy,T(Sy))\leq kd(y,Sy)+(1-k)d(A,B).
 \end{equation}
  If $x$ is a $G$-approximate fixed point for $T$, then there exists a $\epsilon>0$ and by \ref{3.3}
\begin{eqnarray*}
d(Tx, S(Tx)) &\leq & k d(x,Tx)+(1-k)d(A,B)\\
  &\leq& k(d(A,B)+\epsilon)+(1-k)d(A,B)\\
  &=& d(A,B)+k\epsilon \\
  &<& d(A,B)+\epsilon.
  \end{eqnarray*}
  and $(x,Tx)\in P^{Ga}_{(T,S)}(A,B),$ also if $y$ is a $G$-approximate fixed point for $S$, then there exists a $\epsilon>0$ and by \ref{3.4}
  \begin{eqnarray*}
  d(Sy,T(Sy)) &\leq& kd(y,Sy)+(1-k)d(A,B)\\
  &\leq& k(d(A,B)+\epsilon)+(1-k)d(A,B)\\
  &=& d(A,B)+k\epsilon \\
  &<& d(A,B)+\epsilon.
  \end{eqnarray*}
  and $(y,Sy)\in P^{Ga}_{(T,S)}(A,B)$. Therefore $P^{Ga}_{(T,S)}(A,B)\neq \emptyset.$ $\blacksquare$
  
  \begin{theorem} Let $A$ and $B$ be nonempty subsets of a metric space endowed with a graph $G$ and$T:A\cup B \rightarrow A\cup B$, $S:A\cup B \rightarrow A\cup B$ be two continuous maps such that  $T(A)\subseteq B$, $S(B)\subseteq A.$ We suppose that:\\
 i) for every $(x,y)\in E(G),$
$ d(Tx,Sy)\leq \alpha d(x,y)+\gamma d(A,B)$ where $\alpha,\gamma \geq0$ and $\alpha+\gamma=1$;\\
ii)  for any $\{x_n\}_{n\in N}$  and $\{y_n\}_{n\in N}$ in $A\cup B$ as flowing:
 \begin{equation*}
 x_{n+1}=Sy_n~~~,~~~y_{n+1}=Tx_n~~~~for~~some ~~~(x_{1},y_{1})\in E(G)~ for ~n\in N.
 \end{equation*}
 for any $\{x_n\}_{n\in N}$ in $A\cup B$, if $x_n\rightarrow x$ and $(x_n,x_{n+1})\in E(G)$ for $n\in N$ then there is a subsequence convergent $\{x_{n_k}\}_{k\geq 1}$ with $(x_{n_k},x)\in E(G)$ for $n\in N$.\\ Then there exists a $x \in A$ such that $d(x,Tx)=d(A,B).$
\end{theorem}

 {\bf Proof.} We have
\begin{eqnarray*}
  d(x_{n+1},y_{n+1}) &=&  d(Tx_n,Sy_n)\\
  &\leq& \alpha d(x_n,y_n)+\gamma (d(A,B)\\
  &\leq& ...\\
  &\leq& \alpha ^{n+1}d(x_0,y_0)+(1+\alpha +...+\alpha^n )\gamma d(A,B).
\end{eqnarray*}
If $\{x_{n_k}\}_{k\geq 1}$ is converge to $x_1\in A$, that is $ x_{n_k}\rightarrow x_1$ with $(x_{n_k},x_1)\in E(G)$ for all $n \in N$. Then $$d(x_{n_{K+1}},y_{n_{k+1}})\leq \alpha ^{n_{k+1}}d(x_0,y_0)+(1+\alpha +...+\alpha^{n_k} )\gamma d(A,B)$$
Since T is continuous, then $$d(x_{n_{k+1}},Tx_{n_k})\rightarrow \frac{\gamma}{1-\alpha}d(A,B)=d(A,B).$$
Therefore $d(x_1,Tx_1)=d(A,B)$. $\blacksquare$

\begin{definition} Let $A$ and $B$ be nonempty subsets of a metric space endowed with a graph $G$ and Let $T:A\cup B \rightarrow A\cup B$ , $S:A\cup B\rightarrow A\cup B$ be continues maps such that  $T(A)\subseteq B$ , $S(B)\subseteq A$. We define diameter $P^{Ga}_{(T,S)}(A,B)$ by,
$$diam(P^{Ga}_{(T,S)}(A,B))=\sup \{d(x,y):~~d(Tx, Ty)\leq \epsilon+d(A,B)~~for ~some ~~ \epsilon>0\}.$$
\end{definition}

\begin{example} Suppose $A=\{ (x,0):~0\leq x\leq 1\}$, $B=\{(x,1):~0\leq x\leq 1 \}$,
$T(x, 0)=T(x,1)=(\frac{1}{2}, 1)$ and $S(x, 1)=S(x,0)=(\frac{1}{2},0)$. Then $d(T(x, 0), S(y, 1))=1$ and
$diam (P^{Ga}_{(T,S)}(A,B)=diam (A\times B)=\sqrt{2}$.
\end{example}
 \begin{theorem} {\it Let $A$ and $B$ be nonempty subsets of a metric space endowed with a graph $G$ and Let $T:A\cup B \rightarrow A\cup B$ , $S:A\cup B\rightarrow A\cup B$ be continues maps such that $T(A)\subseteq B$ , $S(B)\subseteq A$. If, there exists a $k \in [0, 1],$
$$ d(x, Tx)+d(Sy, y)\leq kd(x, y).$$
Then} $$diam (P^{Ga}_{(T,S)}(A,B) \leq \frac{\epsilon}{1-k}+\frac{d(A,B)}{1-k}~for~some~\epsilon >0.$$
\end{theorem}

{\bf Proof.} If $(x, y) \in P^{Ga}_{(T,S)}(A,B)$, then
\begin{eqnarray*}
d(x, y) &\leq& d(x, Tx)+d(Tx, Sy)+d(Sy, y)\\
  &\leq& \epsilon +k d(x, y)+d(A,B).
\end{eqnarray*}
Therefore $d(x, y) \leq \frac{\epsilon}{1-k}+\frac{d(A,B)}{1-k}.$
Then $diam (P^{Ga}_{(T,S)}(A,B))\leq \frac{\epsilon}{1-k}+\frac{d(A,B)}{1-k}.~\blacksquare $\vspace{0.4cm}

\ConflictofInterests
\section{\bf\large Conclusions}

Nowadays, fixed points, approximate fixed points and graphs in metric spaces play an important role in different areas of mathematics and its applications, particularly in mathematics, physics, biochemistry, electrical engineering and computer science. 
In this work,  we introduced  the new classes of operators and  contraction maps and gave results about $G$-approximate best proximity points and diameter $G$-approximate best proximity points on the metric space $X$ endowed with a graph $G$ such that the set $V (G)$ of vertices of $G$ coincides with $X.$ Also, by using two general Theorem regarding  $G$-approximate best proximity points  of  cyclic maps on metric spaces we  proved several  $G$-approximate best proximity points
 theorems and diameter $G$-approximate best proximity points for a new class of operators and contraction mapping on the metric space $X$ endowed with a graph $G$. We accompanied our theoretical results by some applied examples. 

\end{document}